\documentclass{mpi2015-cscpreprint}

\usepackage[american]{babel}
\usepackage{amssymb}
\usepackage{amsthm,amsmath}
\usepackage{algorithm}%
\usepackage{algorithmicx}%
\usepackage{algpseudocode}%
\usepackage{graphicx}
\usepackage{appendix}
\usepackage{multirow}
\usepackage{array}
\usepackage{amsmath,amssymb}
\usepackage{graphics,graphicx,epsfig,algorithm,algpseudocode}
\usepackage{subcaption}
\usepackage{booktabs}

\usepackage{acronym}
\usepackage{bookmark}
\usepackage{adjustbox}
\usepackage{csquotes}
\usepackage{dsfont}
\usepackage[font=small,labelfont=bf]{caption}
\usepackage{todonotes}
\usepackage{mathtools}
\usepackage{tcolorbox}
\usepackage{hyperref}

\numberwithin{equation}{section}

\def\R{\ensuremath{\mathbb{R}}}
\def\bb{\ensuremath{\mathbf{b}}}
\def\br{\ensuremath{\mathbf{r}}}
\def\bu{\ensuremath{\mathbf{u}}}
\def\bv{\ensuremath{\mathbf{v}}}
\def\bx{\ensuremath{\mathbf{x}}}

\def\loss{\mathcal L}
\def\batch{\mathcal B}

\begin{document}

\title{Data-Driven Modeling of a Controlled Orthotropic Plate Using Machine Learning}

\author[1]{Y. Kim}
\author[3]{A. Zuyev$^{1,}$}
\author[2]{F. Pellicano}
\author[2]{A. Zippo}
\affil[1]{Max Planck Institute for Dynamics of Complex Technical Systems, Magdeburg, Germany
    \email{ykim@mpi-magdeburg.mpg.de}, \orcid{0000-0003-4181-7968}
    \email{zuyev@mpi-magdeburg.mpg.de}, \orcid{0000-0002-7610-5621}
}
\affil[2]{Department of Engineering ''Enzo Ferrari'', University of Modena and Reggio Emilia, Modena, Italy
\email{francesco.pellicano@unimore.it}, \orcid{0000-0003-2465-6584}
\email{antonio.zippo@unimore.it}, \orcid{0000-0001-6206-2619}
}
\affil[3]{Institute of Applied Mathematics and Mechanics, National Academy of Sciences of Ukraine}




\abstract{
We study the problem of learning the input-output map of a controlled vibrating plate with a composite structure from experimental measurements. Analytical modeling of this control system faces challenges due to the essential orthotropy and unknown damping characteristics of the material. Surrogate models based on linear regression, multilayer perceptrons, and gated recurrent units are constructed from the available sampled data. Through comparative analysis, we show that the multilayer perceptron model provides an acceptable approximation of this dynamical system, capturing the potentially nonlinear phenomena in its input-output behavior.
}

 \novelty{
\begin{itemize}
\item Development and experimental comparison of surrogate models for an orthotropic controlled elastic plate, where anisotropic properties, nonlinear actuation characteristics, and complex boundary conditions play a significant role in the dynamics.
\item The comparison of linear regression models with multilayer perceptrons (MLPs) and gated recurrent unit (GRU) techniques shows that nonlinearity plays a significant role in the system’s behavior.
\item It is justified that a moderate artificial neural network depth ($h = 4$ or $6$) is sufficient to effectively leverage the experimental input–output behavior, and that the pretrained MLPs have acceptable approximation ability with the suggested time history length of $s = 200$ samples.
\end{itemize}}

\maketitle

\section{Introduction}\label{sec:intro}
Modeling engineering systems from first principles and using surrogate data-driven approaches represent two different paradigms.
Purely numerical surrogate methods often fail to provide sufficient insight into the influence of mechanical parameters and the effects of time-varying input actions, since they are typically trained as input–output mappings and therefore lack explicit physical structure or interpretability (see, e.g.,~\cite{schar2024emulating,zuyev2023approximating} and references therein).
Although analytic modeling based on variational formulations provides valuable insight into the physics of vibrating flexible structures~\cite{berd,ZS15,zuyev2015partial}, an accurate representation of damping, structural nonlinearities, and actuator–sensor dynamics typically requires additional analysis and often relies on system identification techniques and experimental data~\cite{chopra2020dynamics,inman2013engineering,noel2018grey,al2020critical,lopez2025data}.

\begin{figure}[h]
\centering
\includegraphics[width=0.6\columnwidth]{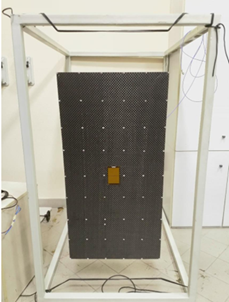}
\caption{Controlled orthotropic plate at the Enzo Ferrari Department of Engineering, University of
Modena and Reggio Emilia.}
\label{fig:plate}
\end{figure}

The goal of this paper is to determine the range of applicability of surrogate modelling techniques for a specific flexible structure -- controlled orthotropic plate  -- and to identify machine learning models that achieve efficient performance and have potential for future applications in control design.
High-fidelity data obtained from an experimental setup at the Enzo Ferrari Department of Engineering, University of Modena and Reggio Emilia, serve as the basis for a comparative analysis of different machine learning approaches. The experimental setup is shown in Fig.~\ref{fig:plate} and described in Section~\ref{sec:setup}. To provide further insight into the mechanical characteristics of the vibrating system and to highlight modelling challenges in the context of an analytical input–output description, we present a previously developed linear PDE-based model of the plate in Section~\ref{sec:rel}.


To address possible nonlinearities in the input-output behavior of the considered mechanical system and to present computationally feasible validation procedures,
we summarize the construction of surrogate models in Section~\ref{sec:model}.
Section~\ref{sec:experi} presents extensive numerical results comparing the input–output behavior of the orthotropic plate using different machine-learning-based surrogate models.
The principal outcomes of our study are summarized in the Conclusion.

\section{Experimental Setup}\label{sec:setup}

The experiments were conducted on a rectangular composite sandwich plate with orthotropic behavior, having in-plane dimensions of $1000 \times 500\,mm$ (see Fig.~\ref{fig:plate}). The plate consisted of two outer carbon-fiber/epoxy laminated skins and a $3\,mm$-thick NOMEX honeycomb polymeric core; the carbon composite was identified as GG285P(T700)-DT120-40, and the fibers were woven in a 0/90 configuration. The specimen was mounted on a dedicated support frame by means of thin nylon wires passing through four holes drilled at the plate corners and connected to elastic bands fixed to the frame members. The suspension geometry and the tensioned lengths of the elastic elements were made as uniform as possible in order to approximate free-boundary conditions and to limit the influence of the supporting structure on the measured dynamics.

The plate was excited by an electrodynamic shaker installed behind the specimen and suspended from the supporting structure in order to reduce parasitic mechanical constraints. The shaker was connected to the plate through a slender stinger so as to provide an approximately pointwise transverse excitation while minimizing additional structural constraints. Similar experimental strategies for vibrating plates and flexible structures have been adopted in previous investigations on orthotropic structures and data-driven modeling approaches. In the present study, particular attention was devoted to the calibration of the actuator–sensor chain in order to ensure consistency between the applied input and the measured response. The excitation force generated by the shaker was monitored by a load cell placed along the excitation path, between the shaker rod and the plate surface, allowing the input signal $u(t)$ to be reconstructed from the measured voltage output of the force sensor after appropriate signal conditioning, amplification, scaling, and calibration. The acquisition channels were organized so that the load-cell signal represented the input and the piezoelectric response represented the output.

The structural response of the plate was measured through a bonded piezoelectric transducer used as a strain gauge, whose voltage output was taken as the experimental observable. The sensor was installed at a non-collocated location with respect to the actuator, thus providing a sensitive measurement of the local curvature of the plate and enabling the detection of bending vibrations associated with the orthotropic elastic properties of the composite structure. The setup documentation identifies the transducer as a Macro Fiber Composite manufactured by Smart Materials (MFC P1 type, model M8557-P1), with an active area of $85 \times 57\, mm$, overall dimensions of $103 \times 64\, mm$, capacitance of $12.84\, nF$, free strain of $1350\, ppm$, and blocking force of $693\, N$. Owing to its very small thickness, the transducer could be surface-bonded to the plate using epoxidic resin with negligible modification of the structural properties, while still ensuring sensitive measurement of the local structural deformation.

The actuation and measurement chain was integrated with a digital acquisition and control architecture based on \texttt{dSPACE} and \texttt{MATLAB/Simulink}, ensuring accurate synchronization between excitation and response measurements. The sampling frequency was set to $30\, kHz$ in order to capture the relevant dynamic content of the plate vibrations and avoid aliasing effects in the recorded time histories. In the broader setup, the piezoelectric transducer is driven through a TREK PA05039 high-voltage amplifier, which converts an input in the range from $-2.5\, V$ to $7.5\, V$ into an output from $-500\, V$ to $1500\, V$, with nominal input impedance of $25 \,k\Omega$, output impedance below $0.1\, k\Omega$, and output noise below $50 \,mV$ rms. The electromagnetic shaker was powered by a TIRA BAA 60 amplifier, capable of supplying up to $2.5\, A$ with a maximum $5\, V$ input signal.

During the experiments, the structure was excited by short voltage pulses applied to the shaker at regular time intervals. This excitation strategy generated transient responses followed by freely decaying vibrations governed by the intrinsic structural and material damping of the composite plate. Such a procedure provides a rich dataset containing both the forced response and the subsequent free-vibration regime, which is particularly suitable for modal characterization, dynamical model identification, and the training of surrogate machine-learning models. The resulting dataset therefore contains high-resolution synchronized time histories of the input force and the measured piezoelectric strain response, forming the basis for the data-driven modeling and validation procedures described in the following sections.

\section{Governing Equations of the Orthotropic Plate}\label{sec:rel}

In the paper~\cite{zuyev2024dynamic}, a mathematical model of a controlled orthotropic elastic plate
is studied based on the following Kirchhoff-type partial differential equation:
\begin{equation}\label{Kirchhoff}
\begin{aligned}
&\frac{\partial^2 w(x,t)}{\partial t^2} + \alpha \frac{\partial w(x,t)}{\partial t} + d_{11} \frac{\partial^4 w(x,t)}{\partial x_1^4}    +d_{22} \frac{\partial^4 w(x,t)}{\partial x_2^4} \\
&+ 2(d_{12}+2 d_{66})\frac{\partial^4 w(x,t)}{\partial x_1^2 \partial x_2^2} =  \frac{\delta_{S_0}(x)}{\rho h} u(t),\; x\in D,\;t\ge 0.
\end{aligned}
\end{equation}
The plate is assumed to be rectangular, and its middle surface is parameterized by the spatial coordinates $x=(x_1,x_2)\in D=[0,\ell_1]\in [0,\ell_2]$,
where $\ell_1$ and $\ell_2$ denote the side lengths of the rectangle $D$.
The function $w=w(x,t)$ describes the transverse displacement of the plate's middle surface at time~$t$.

The plate is controlled by the force $u(t)$ applied by an actuator at the point $S_0\in D$,
and its action is represented by the right-hand side of~\eqref{Kirchhoff} through the Dirac delta function concentrated at $S_0$.
The state of this plate is measured by a piezoelectric sensor at the point $\bar x \in D$, which provides the output signal of the form
\begin{equation}\label{output}
y(t) = k w_{x_1 x_1} (t,\bar x),
\end{equation}
where $k>0$ corresponds to the sensor gain.

The Kirchhoff equation~\eqref{Kirchhoff} is considered together with the following boundary conditions at $x_1=0$ and $x_1=\ell_1$:
\begin{equation}\label{BC_F}
\begin{aligned}
& w_{x_1 x_1} + \nu_2 w_{x_2 x_2} = 0, \\
& w_{x_1 x_1 x_1}+\left( \nu_2 + \frac{4(1-\nu_1\nu_2)G}{E_1} \right) w_{x_1 x_2 x_2} =0,
\end{aligned}
\end{equation}
and the boundary conditions at
$x_2=0$ and $x_2=\ell_2$:
\begin{equation}\label{BC_x2}
\begin{aligned}
&\nu_1 w_{x_1 x_1} + w_{x_2 x_2} = 0, \\
& \left( \nu_1 + \frac{4(1-\nu_1\nu_2)G}{E_2} \right) w_{x_1 x_1 x_2} + w_{x_2 x_2 x_2} =0.
\end{aligned}
\end{equation}
These conditions correspond to free boundary conditions, i.e., the absence of boundary forces and bending moments.

The coefficients of the control system~\eqref{Kirchhoff}--\eqref{BC_x2} are discussed in~\cite{zuyev2024dynamic},
where mechanically consistent values of the parameters, encoding the real material properties of the experimental setup shown in Fig.~\ref{fig:plate}, are reported as follows:
$$
\begin{aligned}
& \ell_1 = 1\, m,\; \ell_2= 0.5\, m,\, h=3.6\,mm,\, \rho = 505.6\,\frac{kg}{m^3},\\
& d_{11} = \frac{E_1 h^2}{12\rho(1-\nu_1 \nu_2)}, \; d_{22}= \frac{E_2 h^2}{12\rho(1-\nu_1 \nu_2)}, \\\
& d_{12} = \nu_2 d_{11}=\nu_1 d_{22},\;  d_{66} = \frac{G h^2}{12 \rho},\\
& E_1 = 23\,GPa,\;E_2 = 14\,GPa,\, G=2.2\,GPa,\\
& \nu_1 = 0.25,\,\nu_2 = \frac{\nu_1 E_2}{E_1},\\
& S_0 = (0.17 \,m,0.25 \, m),\; C_0 = (0.5 \,m,0.21\,m).
\end{aligned}
$$

For the linear PDE-governed control system~\eqref{Kirchhoff}--\eqref{BC_x2}, a family of finite-dimensional approximations is constructed in~\cite{zuyev2024dynamic} using Galerkin's method.
The frequency-response characteristics, represented by the transfer function, are derived in the aforementioned paper, and the basic properties of the transfer functions are compared between the analytical and data-driven models. It should be noted that the damping coefficient $\alpha$ in~\eqref{Kirchhoff}, corresponding to a viscous damping model, is not known a priori; therefore, adequately modeling the damping requires further analysis.

As can be seen, the control system described by the fourth-order partial differential equation~\eqref{Kirchhoff}, boundary conditions involving third-order derivatives~\eqref{BC_F}--\eqref{BC_x2}, input dependence through the Dirac delta function, and the output map~\eqref{output} defined in terms of a second-order partial derivative already constitutes a mathematically complex model containing multiple parameters. The complexity of the model description remains significant even in the Galerkin approximations, partly due to the presence of damping and gain parameters that must be identified.

It should also be noted that the above analytical models are linear and therefore cannot capture nonlinearities if they are present in the measured input–output response. For this reason, in Section~\ref{sec:model} we explore an alternative approach based on surrogate models. Such models can potentially capture nonlinear behavior and allow us to obtain computationally feasible results for model validation using the available experimental data.

\section{Surrogate Modeling}\label{sec:model}
Within the considered surrogate modelling framework, we assume that
the output signal $y(t)\in \mathbb R$ at time $t\ge 0$ is affected not only by the current input $u(t)\in \mathbb R$, but also by a sequence of preceding inputs, since the plate's oscillation is governed by the cumulative influence of its input history.
To capture this behavior in surrogate modeling, we construct time series surrogate models represented as
$$
\hat y(t)=f(u(t), u(t-\Delta t),\cdots , u(t-(s-1)\Delta t);\theta),
$$
where $\hat y(t)\in\mathbb R$ denotes the predicted system output,
$\Delta t>0$ is the time step, $s\in{\mathbb N}$ denotes the sequence length,
and the vector $\theta\in{\mathbb R}^p$ comprises all parameters of $f$ to be identified.
Equivalently, we express the model as
\begin{equation}\label{sur_mod}
\hat y(t)=f({\bu}_s(t);\theta)
\end{equation}
where
\begin{equation}\label{us}
{\bu}_s(t)=(u(t), u(t-\Delta t),\cdots , u(t-(s-1)\Delta t))^\top\in\R^s
\end{equation}
The function $f:{\mathbb R}^s\times {\mathbb R}^p\to \mathbb R$ and its parameters $\theta\in\mathbb R^p$ are to be identified from the available experimental data.

As a starting point, we consider linear regression models (LRs) described by
\begin{equation}\label{LR}
\hat y(t) = W{\bu}_s(t) + b,
\end{equation}
where $W\in\R^{1\times s}$ is a weight matrix and $b\in\R$ is a bias term.
This model provides a baseline for capturing input–output dependencies and serves as a foundation for further exploration of nonlinear modeling approaches.

Next, we develop multilayer perceptrons (MLPs) consisting of fully connected layers, each followed by a nonlinear activation function in the hidden layers. An MLP is a fundamental architecture of feedforward neural networks that constructs a nonlinear mapping from the input space to the output space through successive compositions of affine transformations interleaved with nonlinear activation functions~\cite{rumelhart1987learning,Goodfellow-et-al-2016}. The inclusion of nonlinearity allows the model to effectively capture nonlinear relationships between the model input ${\bu}_s(t)$ and the output $\hat y(t)$, thus potentially providing higher accuracy than linear approaches.

At each $t\ge s\Delta t$, the input-output map modelled by the MLPs is presented as
\begin{equation}\label{MLP}
\begin{aligned}
\bx_1(t) &= \text{ELU}(W_1\bu_s(t)+\bb_1), \\
\bx_2(t) &= \text{ELU}(W_2\bx_1(t)+\bb_2), \\
\vdots & \\
\bx_h(t) &= \text{ELU}(W_h\bx_{h-1}(t)+\bb_h), \\
\hat y(t) &= W\bx_h(t)+b,
\end{aligned}
\end{equation}
where $W_i$ and $\bb_i$ are
the weight matrix and bias vector in the
$i$-th hidden layer, respectively, $i\in\{1,2, ..., h\}$.
The matrix $W$ and the scalar $b$ denote the weight matrix and bias in the output layer.
As can be seen from the above formulas, the $s$-dimensional vectors $\bx_1$, ...,  $\bx_{h}$ are computed at each $t$,
whereas the constant parameters $W_i$, $\bb_i$, and $b$ are determined during the training of the model.
In this paper, we use the Exponential Linear Unit activation function $\textrm{ELU}:{\mathbb R}^s\to{\mathbb R}^s$ which acts
componentwise on a vector $\bx = (x_1,...,x_s)^\top$ as~\cite{Djo16ELU}:
$$
\textrm{ELU}_j(\bx) = \begin{cases}x_j,\; x_j>0,\\ e^{x_j}-1,\; x_j\le 0.\end{cases}
$$

As an alternative to MLPs, various types of recurrent neural networks designed for sequential data can be considered.
In this work, we utilize gated recurrent units (GRUs)~\cite{Ch14GRU},
where the input sequence is processed sequentially by a stack of $h$ layers.
For a given $t\ge s\Delta t$, the state of the zero-th layer is defined as
\begin{equation}\label{GRU_input}
\bx^{(0)}_j(t) = {u_s}_j(t),\quad j=1,..., s.
\end{equation}
Then, for each layer $i\in \{1,\ldots,h\}$ and time step $j\in\{1,\ldots,s\}$, the GRU computational scheme is defined as
\begin{equation}\label{GRU}
\begin{aligned}
\bv_j^{(i)} &= \sigma(W_v^{(i)}\bx_{j-1}^{(i-1)}+V_v^{(i)} \bx_{j}^{(i-1)}+\bb_v^{(i)}), \\
\br_j^{(i)} &= \sigma(W_r^{(i)}\bx_{j-1}^{(i-1)}+V_r^{(i)} \bx_{j}^{(i-1)}+\bb_r^{(i)}), \\
\tilde \bx_j^{(i)} &= \tanh(W_x^{(i)}(\bx_{j-1}^{(i-1)}\odot\br_j^{(i)})+V_x^{(i)} \bx_{j}^{(i-1)}+\bb_x^{(i)}),  \\
\bx_j^{(i)} &= \bv_j^{(i)}\odot\bx_{j-1}^{(i)} +  (1-\bv_j^{(i)})\odot\tilde\bx_j^{(i)},
\end{aligned}
\end{equation}
where $W_v^{(i)}$, $V_v^{(i)}$, $W_r^{(i)}$, $V_r^{(i)}$, $W_x^{(i)}$, and $V_x^{(i)}$ are weight matrices, and $\bv_j^{(i)}$ and $\br_j^{(i)}$ denote the outputs of the update gate and the reset gate, respectively, and $\bx_j^{(i)}\in \mathbb R^{n_i}$ is the hidden state at time step $j$ in the $i$-th hidden layer, $j\in\{1,2, ... , s\}$, $i\in\{1,2, ... , h\}$.
Note that the dimension of the $i$-th layer is $n_i\in \mathbb N$;
so, we emphasize vector quantities using bold notation here and recall that $n_0=1$ for the considered SISO system.
In the above formulas, $\odot$ denotes elementwise multiplication, and $\sigma$ is the sigmoid function, whose componentwise action is given by
$$
\sigma(x)= \frac{1}{1+e^{-x}}\quad \text{for}\;\;x\in{\mathbb R}.
$$

After computing the last hidden state $\bx^{(h)}_s$ in~\eqref{GRU}, the output of the system at time $t$ is represented as
\begin{equation}\label{GRU_output}
\hat y(t) = W \bx_{s}^{(h)} + b.
\end{equation}

The above-defined GRU architecture introduces gating mechanisms that regulate the flow of temporal information across finite time steps.
The update gate $\bv_j^{(i)}$ determines the extent to which the previous hidden state $\bx_{j-1}^{(i)}$ is retained, thus adjusting long-term dependencies.
The reset gate $\br_j^{(i)}$ modulates the influence of past information when computing the candidate hidden state $\tilde{\bx}_j^{(i)}$.

\section{Numerical Experiments}\label{sec:experi}
\begin{figure}[t]
\centering
\includegraphics[width=1\columnwidth]{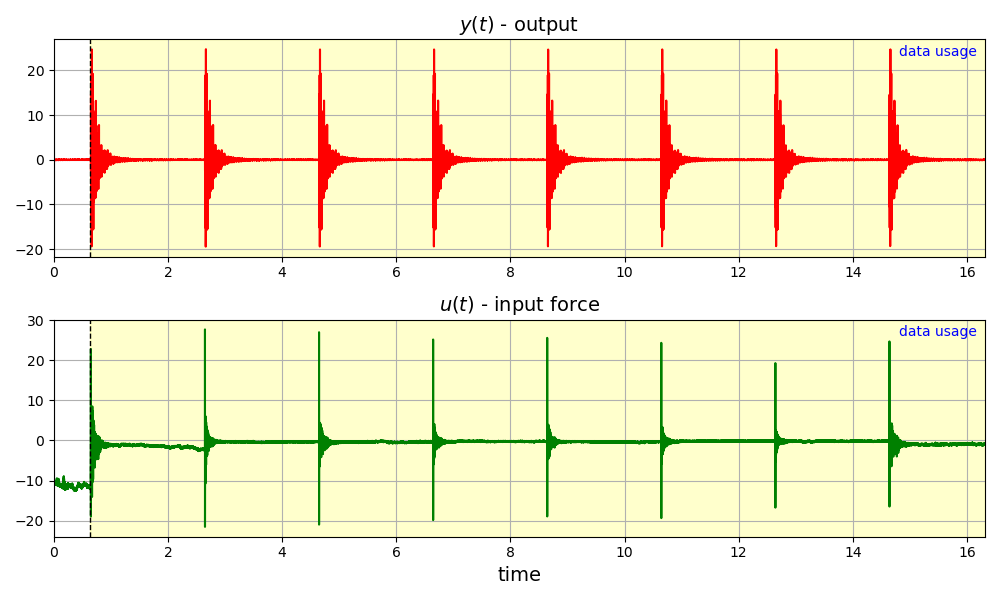}
\caption{Time series dataset with input $u(t)$ and output $y(t)$, sampled over time $t$. Only the data within the yellow-shaded region ($t\in[0.64\,s,16.31\,s]$) are used for training, validation, and testing models.}\label{fig:data}
\end{figure}

We apply the surrogate modeling approaches, described in Section~\ref{sec:model}, to experimental
input-output data obtained from the controlled orthotropic plate setup with a non-collocated sensor-actuator pair.
In this setup, the plate is actuated by an electromagnetic shaker via a voltage input,
which can be physically controlled using \texttt{dSPACE} with digital communication with \texttt{MATLAB/Simulink}.

To measure the force of the actuation, a load sensor is placed between the shaker rod and the plate.
The voltage output of this sensor provides the control signal $u(t)$ of our force-controlled model.

On the opposite side of the plate, a piezoelectric strain gauge is installed such that, with appropriate amplification in the measurement chain, the strain is measured as the output signal $y(t)$.
This measurement scheme is implemented using \texttt{MATLAB} and \texttt{Simulink} with the signals $u(t)$ and $y(t)$
(up to scaling constants matching the amplifiers and ADC/DAC characteristics) recorded at a sampling frequency of $30\, kHz$, corresponding to a sampling interval
 $\Delta t = 0.03(3)\, ms$.

In the conducted experiments, the shaker is actuated by short voltage pulses of duration $\Delta t$, repeated every $2\,s$.
The overall dataset consists of $N_s=489,562$ samples of the form $(t_k,u (t_k),y (t_k))$ with $t_k=k \Delta t$  for $k=0,1,..., N_s-1$.
The time plots of the input and output signals are sketched in Fig.~\ref{fig:data}.
We exclude from the analysis the portion of the time series before the first control impulse at $T_0\approx 0.64\,s$, as shown in Fig.~\ref{fig:data}.

The remaining dataset comprises $470,354$ samples $(t_k, u(t_k), y(t_k))$ in the time domain $t\in [0.64\,s,16.31
\,s]$, which are divided into $316,078$ data points (64\%) in $t\in [0.64\,s,11.17\,s)$ for training, $60,205$ data points (16\%) in $t\in [11.17\,s,14.31\,s)$ for validation, and $94,071$ data points (20\%) in $t\in [14.31\,s,16.31\, s]$ for testing.
These datasets are mutually exclusive and contain no temporal overlap.

As a standard supervised learning approach for training MLPs and GRUs,
we employ the mean squared loss function
\begin{equation}\label{SS}
\loss(\theta) = \frac{1}{|\batch|}\sum_{(t_k,\bu_s(t_k),y_k)\in\batch} ( y_k-f(\bu_s(t_k);\theta) )^2
\end{equation}
where $\batch$ denotes a batch drawn from the dataset and $|\batch|$ represents the batch size, i.e.,
the number of data samples $(t_k,\bu_s(t_k),y_k)$ in  $\batch$, such that the vector $\bu_s(t_k)$ of dimension $s$ is constructed from the dataset by the rule~\eqref{us}.
For each particular LR, MLP, or GRU model, we follow the convention that the parameter vector $\theta$ in~\eqref{SS} and~\eqref{sur_mod} comprises all parameters of the considered model.

To avoid overfitting during neural network training, we adopt a batch training strategy using the \texttt{Adam} optimizer \cite{Kin15ADAM} combined with early stopping based on the validation loss.
In this approach, the training process is terminated if the validation loss does not improve for a predefined number of epochs. The batch size is set to 256 in all experiments.

To evaluate model accuracy, we use the $R^2$ score:
$$R^2 = 1-\frac{\sum_{i=1}^N(y_i-\hat y_i)^2}{\sum_{i=1}^N(y_i-\bar y)^2},$$
where $N$ denotes the number of considered samples, $y_i$ is the reference output, $\hat y_i$ is the predicted output corresponding to $y_i$, and $\bar y=\frac{1}{N}\sum_{i=1}^N y_i$ is the mean of the reference outputs $y_1,y_2, ... , y_N$.
For convenience of notation, we index these data from $1$, as this can always be achieved by an appropriate shift of time without loss of generality.
The methods presented here are implemented in \texttt{Python} and are available in the GitHub repository~\cite{py_code}.
%
%


\begin{table*}[t]
\centering
\makebox[\textwidth][c]{}
\begin{subtable}[t]{\textwidth}
\centering
\begin{tabular}{@{}l l@{}}
(a) \textbf{LR}$\quad$ &
\begin{tabular}{c|cccc}
\toprule
$h$ \textbackslash~$s$ & 10 & 50 & 100 & 200 \\
\midrule
0 & NS / NS & 0.016 / 0.011 & 0.039 / 0.027 & \textbf{0.089 / 0.111} \\
\bottomrule
\end{tabular}
\end{tabular}
\end{subtable}
\hspace{0.03\textwidth} 
\begin{subtable}[t]{\textwidth}
\centering
\begin{tabular}{@{}l l@{}}
(b) \textbf{MLP} &
\begin{tabular}{c|cccc}
\toprule
$h$ \textbackslash~$s$ & 10 & 50 & 100 & 200 \\
\midrule
1 & 0.017 / 0.006 & 0.404 / 0.271 & 0.594 / 0.416 & 0.774 / 0.533 \\
2 & 0.042 / 0.023 & 0.787 / 0.563 & 0.891 / 0.778 & 0.941 / 0.829 \\
4 & 0.058 / 0.039 & 0.861 / 0.692 & 0.945 / 0.847 & 0.975 / 0.909 \\
6 & 0.061 / 0.007 & 0.868 / 0.675 & 0.952 / 0.842 & \textbf{0.977 / 0.908} \\
\bottomrule
\end{tabular}
\end{tabular}
\end{subtable}
\hspace{0.03\textwidth} 
\begin{subtable}[t]{\textwidth}
\centering
\begin{tabular}{@{}l l@{}}
(c) \textbf{GRU} &
\begin{tabular}{c|cccc}
\toprule
$h$ \textbackslash~$s$ & 10 & 50 & 100 & 200 \\
\midrule
1 & 0.054 / 0.015 & 0.719 / 0.404 & 0.778 / 0.571 & 0.888 / 0.732 \\
2 & 0.085 / 0.018 & 0.752 / 0.417 & 0.813 / 0.548 & \textbf{0.913 / 0.814} \\
4 & 0.092 / 0.016 & 0.502 / 0.204 & 0.804 / 0.547 & 0.774 / 0.495 \\
6 & 0.034 / 0.027 & 0.260 / 0.168 & 0.005 / 0.005 & 0.454 / 0.322 \\
\bottomrule
\end{tabular}
\end{tabular}
\end{subtable}
\caption{$R^2$ score (training/test): $s$ denotes the input sequence length, and $h$ is the number of hidden layers. The results correspond to the best performance of each model,
determined by the highest training $R^2$ score across 3 runs.}\label{tab:results}
\end{table*}

Table~\ref{tab:results} presents the training and test $R^2$ scores for LR~\eqref{LR}, MLP~\eqref{MLP}, and GRU~\eqref{GRU_input}--\eqref{GRU_output} models across various input sequence lengths $s\in\{10,50,100,200\}$ and numbers of hidden layers $h$.

As we observe, the LR model~\eqref{LR} performs poorly overall, indicating limited predictive capability, and achieves its best test score of only 0.111 at $s=200$.

In contrast, the MLP model~\eqref{MLP} shows strong performance that improves consistently with both increasing sequence length and network depth.
The best overall result is achieved by the MLP with six hidden layers at $s=200$, attaining training and test $R^2$ scores of 0.977 and 0.908 respectively.

The GRU model~\eqref{GRU_input}--\eqref{GRU_output} also shows strong performance even with shorter sequences, outperforming both LR and shallow MLPs for
$s=10$. However, its performance does not always improve with depth; deeper GRU architectures tend to suffer from degraded generalization. The GRU achieves its highest test $R^2$ score of 0.814 with two hidden layers at $s=200$.
These results demonstrate that while recurrent models such as GRU are effective at predicting outputs with shallow architectures, their generalization tends to degrade as the depth increases. In contrast, deep feedforward models such as MLP benefit significantly from additional layers when longer input histories are available.

\begin{figure}[t]
\centering
\includegraphics[width=1\columnwidth]{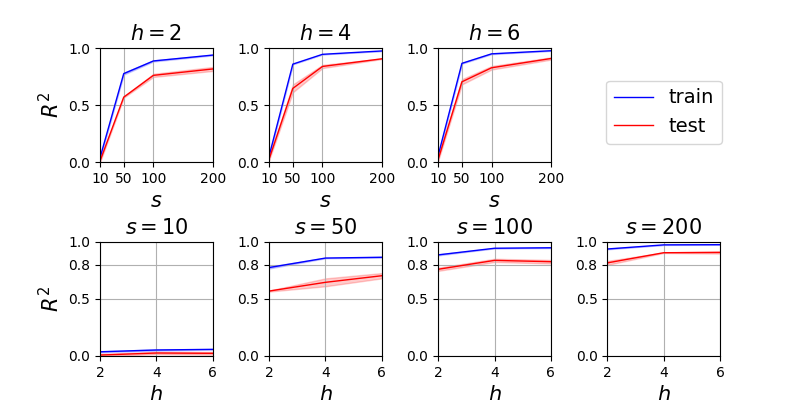}
\caption{$R^2$ scores of pretrained MLPs~\eqref{MLP}: Solid lines represent the mean score $\mu$, and the shaded regions indicate the statistical uncertainty measured using the standard deviation $\sigma$, i.e., the interval $(\mu-\sigma,\mu+\sigma)$.}
\label{fig:mlp-res}
\end{figure}

Fig.~\ref{fig:mlp-res} illustrates the effect of varying the input sequence length $s$ and the number of hidden layers $h$ on the predictive performance of pretrained MLPs, evaluated using the $R^2$ score.
The top row presents $R^2$ scores as a function of $s\in\{10,50,100,200\}$ for fixed hidden layer depths $h=2,4,6$.
Conversely, the bottom row shows performance as a function of the number of hidden layers $h\in{2,4,6}$ for fixed sequence lengths. In each plot, solid lines represent the mean $R^2$ score $\mu$ computed over three independent runs, while the shaded regions depict statistical uncertainty, quantified as one standard deviation $\sigma$, i.e., the interval $\mu\pm\sigma$.
Deeper MLPs generally lead to improved accuracy, particularly at larger $s$, with the best performance achieved at $h=6$ and $s=200$.
However, the marginal benefit of adding layers diminishes beyond a certain point, as shown by the relatively flat curves in the bottom row for higher $s$ values.
Furthermore, the results show a consistent increase in both training and test $R^2$ scores as the input sequence length $s$ increases. In particular, a substantial improvement is observed between $s=10$ and
$s=50$, after which the performance gains begin to plateau.

 \begin{figure*}[t]
    \centering
    \begin{subfigure}[t]{0.32\linewidth} 
        \centering
        \includegraphics[width=\linewidth]{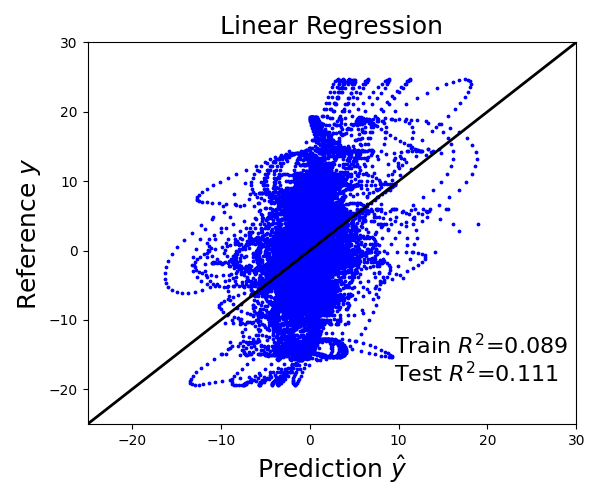}
        \caption{LR ($s=200, h=0$)}
    \end{subfigure}
    \begin{subfigure}[t]{0.32\linewidth}
        \centering
        \includegraphics[width=\linewidth]{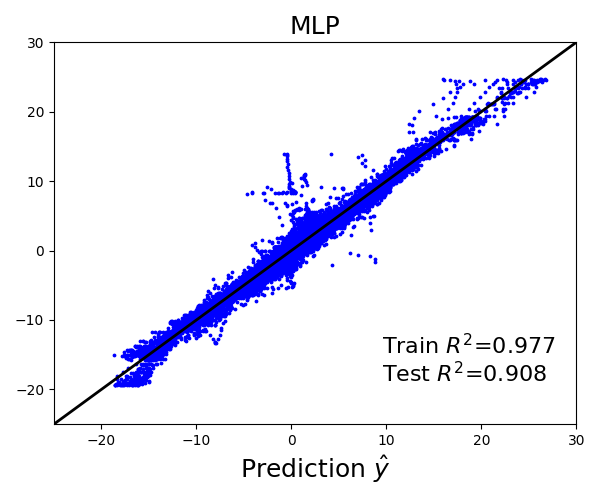}
        \caption{MLP ($s=200, h=6$)}
    \end{subfigure}
    \begin{subfigure}[t]{0.32\linewidth}
        \centering
        \includegraphics[width=\linewidth]{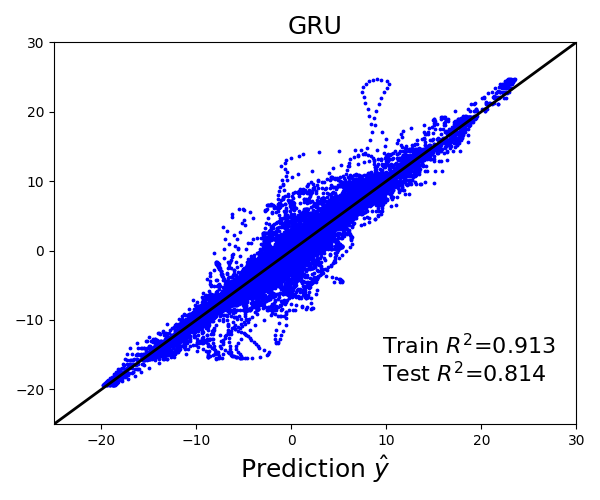} 
        \caption{GRU ($s=200, h=2$)}
     \end{subfigure}
\caption{Scatter plots of predicted versus reference outputs for the LR~\eqref{LR}, MLP~\eqref{MLP}, and GRU~\eqref{GRU_input}--\eqref{GRU_output} models.}
\label{fig:scatter}
\end{figure*}

Fig.~\ref{fig:scatter} illustrates the performance of each model using its best pretrained model weights. The diagonal black line represents the ideal case where $y=\hat y$.
The LR model shows a broad and dispersed distribution of prediction errors, with significant deviation from the identity line, resulting in poor predictive accuracy.
In contrast, the MLP with six hidden layers yields predictions that align closely with the identity line, indicating high accuracy and low variance
for both the training and test data.
The GRU model~\eqref{GRU_input}--\eqref{GRU_output} with two hidden layers also demonstrates improved prediction quality relative to LR~\eqref{LR},
although its predictions exhibit slightly greater spread compared to those of the MLP~\eqref{MLP}.

\section{Conclusion}\label{sec:concl}
To the best of our knowledge, this paper is the first to address the development and experimental comparison of concrete machine learning surrogate models for orthotropic controlled elastic plates, where anisotropic properties, nonlinear actuation characteristics, and complex boundary conditions play a significant role in the dynamics.

Our computational findings, based on the developed code~\cite{py_code}, indicate that increasing the input sequence length significantly enhances model accuracy, and that moderate network depth (e.g., $h=4$ or $6$) is sufficient to leverage this information effectively. Furthermore, the narrow shaded regions in Fig.~\ref{fig:mlp-res} indicate stable and reliable performance for the pretrained MLPs.
A clear indication of approximation efficiency for the time history length $s=200$
 suggests verifiable bounds for possible approximations of the theoretically infinite-dimensional control system using finite-dimensional difference equations. Therefore, the contribution of this work is not limited to surrogate modeling alone.

 Furthermore, the comparison of linear regression models with nonlinear perceptrons and gated recurrent unit techniques shows that nonlinearity plays a significant role in the system’s behavior. Consequently, the linear models exhibit limitations, particularly regarding their range of applicability. Physically, these conclusions are underpinned by the properties of the composite material—which are entirely governed by linear elasticity theory—while the nonlinear characteristics of the actuator-sensor pair are accounted for in the proposed surrogate models.

In future work, we plan to apply the developed modeling techniques to control design for active vibration damping and the development of data-driven optimal control methods.


\section*{Acknowledgments}%
\vskip-1ex
Yongho Kim was supported by the German Research Foundation (DFG) through the research training group 2297 ``MathCoRe", Magdeburg.
Alexander Zuyev gratefully acknowledges the funding by the
European Regional Development Fund (ERDF) within the programme Research
and Innovation — Grant Number ZS/2023/12/182138.

\end{document}